\documentclass[fleqn,11pt]{article}
\usepackage{graphicx}% Include figure files
\usepackage{amsfonts}

\newcommand{\bq}{\begin{equation}}            
\newcommand{\eq}{\end{equation}}            
\newcommand{\bqa}{\begin{eqnarray}}            
\newcommand{\eqa}{\end{eqnarray}}            
\newcommand{\ba}{\begin{array}} 
\newcommand{\ea}{\end{array}} 
\newcommand{\bi}{\begin{itemize}} 
\newcommand{\ei}{\end{itemize}} 

\newcommand{\al}{\alpha}                
\newcommand{\be}{\beta}

\newcommand{\nb}{{\nabla}}
\newcommand{\vs}{\vspace{2mm}}
\newcommand{\hs}{\hspace{5mm}}
\newcommand{\ra}{\rightarrow}           

\newcommand{\df}{\stackrel{\triangle}{=}}     

\begin{document}
%\renewcommand{\thepage}{\roman{page}}
%\pagenumbering{roman}
%\baselineskip=0.65cm
%\baselineskip=0.6cm
\topmargin=-0.5in

\begin{center}
{\Large \bf Stiffness Mitigation in Stochastic Particle Flow Filters}\\[2mm]
\end{center}
\begin{center}
\begin{tabular}{c  c}
        Liyi Dai & Frederick E. Daum\\
        Raytheon Missiles \& Defense & \hspace{10mm} Raytheon Missiles \& Defense \\
        50 Apple Hill Drive & 235 Presidential Way \\
        Tewksbury, MA 01876 & Woburn, MA 01801 \\
        liyi.dai@raytheon.com & daum@raytheon.com 
\end{tabular}
\end{center}
\vs

\begin{abstract}
The linear convex log-homotopy has been used in the derivation of particle flow filters. One natural question is whether it is beneficial to consider other forms of homotopy. We revisit this question by considering a general linear form of log-homotopy for which we derive particle flow filters, validate the distribution of flows, and obtain conditions for the stability of particle flows. We then formulate the problem of stiffness mitigation as an optimal control problem by minimizing the condition number of the Hessian matrix of the posterior density function. The optimal homotopy can be efficiently obtained by solving a one-dimensional second order two-point boundary value problem. Compared with traditional matrix analysis based approaches to condition number improvements such as scaling, this novel approach explicitly exploits the special structure of the stochastic differential equations in particle flow filters. The effectiveness of the proposed approach is demonstrated by a numerical example.
\end{abstract}

{\bf Keywords.} particle flow filters, nonlinear filtering, homotopy, stiffness mitigation, stability, optimal control, two-point boundary value problem

\vspace{0.5cm}

\centerline{\today}
%\centerline{March 22, 2021}

%\newpage

\section{Introduction}
Consider two probability density functions $p_0(x)$ and $p_1(x)$. For each $\lambda \in [0, 1]$, define
\bq
p(x,\lambda) = \frac{p_0^{1-\lambda}(x)p_1^{\lambda}(x)}{\Gamma(\lambda)},
\label{line_homotopy}
\eq
in which $\Gamma(\lambda)>0$ is the normalization factor so that $p(x,\lambda)$ is a probability density function. Note that
\[    p(x,0) = p_0(x), \hs p(x,1)=p_1(x). \]
Therefore, the mapping $p(x,\lambda): \mathbb{R}^+\times [0, 1]\longrightarrow \mathbb{R}^+$ defines a homotopy from $p_0(x)$ to $p_1(x)$. By taking natural logarithm on both sides of (\ref{line_homotopy}), we obtain
\bq
\log p(x,\lambda)=(1-\lambda)\log p_0(x)+\lambda \log p_1(x)-\log \Gamma(\lambda).
\label{line_loghomo}
\eq
In other words, $\log p(x,\lambda)$ is a linear convex combination of $\log p_0(x)$ and $\log p_1(x)$. The straight line log-homotopy (\ref{line_loghomo}) has been used in the derivation of particle flow filters \cite{DH2007,DH2013,DH2015,DHN2010,DHN2016,DHN2018,Khan}. This straight line log-homotopy is easy to implement, leads to computationally efficient particle flow filters, and avoids mathematical singularities as those in Moser's flow in optimal transport. However, it has been recognized that other forms of homotopy is possible \cite{DH2016}. One natural question is: Is it beneficial to consider other forms of homotopy? Several opening problems regarding the choice of homotopy were discussed in details in \cite{DH2016}, e.g., the Open Problems 5, 7, 8, 10, 11. In this paper, we revisit this question and derive particle flow filters for more general forms of homotopy. As demonstrated in the literature on particle flow filters and in transport theory, there is much freedom in the design of particle flow filters. We limit our discussion to a general linear combination of $\log p_0(x)$ and $\log p_1(x)$. Our results suggest that the general form of homotopy could be useful for addressing numerical stability such as reducing the stiffness of the underlying stochastic differential equation that governs particle flows.

Stiff ordinary or stochastic differential equations pose challenges to obtaining numerical approximations of exact theoretical solutions with guaranteed accuracy \cite{CH,HW,Lam}. While it has not been possible to define stiffness in precise mathematical terms in a satisfactory manner, \cite{Lam} refers stiffness as a general phenomenon in which a numerical method is "forced to use a step-length that is excessively small in relation to the smoothness of the exact solution". Stiffness mitigation in particle flow filters has been discussed in \cite{Cro,DH2014,MDD}, with emphasis on improving numerical implementation without changing the straight line log-homotopy.
In designing particle flow filters, the choice of homotopy affects the form of the underlying stochastic differential equation that we need to solve numerically to construct state estimates and confidence intervals. Therefore, the freedom in the choice of homotopy may potentially be used to mitigate the stiffness of the stochastic differential equation. By exploiting the special structure of the stochastic differential equations for the particle flows, we formulate the problem of stiffness mitigation as an optimal control problem. It turns out that the optimal linear log-homotopy can be obtained by solving a simple one-dimensional second order  two-point boundary value problem. Numerical example demonstrates that the proposed approach indeed reduces the stiffness of particle flow equations, which in turn reduces estimation errors and covariance matrix of particle flow filters. 

On the one hand, we would like to derive an optimal homotopy for stiffness mitigation of flow equations. On the other hand, we also need to maintain the stability of particle flows for numerical stability. The stochastic differential equation of particle flows is time-varying, for which eigenvalue-based stability analysis is generally not applicable \cite{Kha}. In \cite{DD2021}, we proposed an approach of Lyapunov stability to the stability analysis of particle flows. In this paper, we extend this approach to the stability analysis of particle flows under general linear log-homotopy. We find conditions to ensure both the stability of particle flows and stiffness mitigation.

The rest of the paper is organized as the follows. In Section 2, we consider a general linear form of log-homotopy. We first derive particle flow filters, and then prove that the particle flow indeed has the desired posterior density function, thus establishing the correctness of the distribution of particle flows. We further obtain conditions under which particle flows are stable. In Section 3, we formulate stiffness mitigation as an optimal control problem by minimizing the condition number of the Hessian matrix of the posterior density function, and obtain optimal linear log-homotopy. In Section 4, we use an example to illustrate the effectiveness of our proposed approach to stiffness mitigation and demonstrate that such stiffness mitigation indeed leads to error reduction in state estimation. Finally, we summarize our findings in concluding remarks in Section 5. For clarity in discussion, all proofs are moved to the Appendix.
\vs

\noindent
{\em Notations}.

We use $\mathbb{R}^n$ to denote the real-valued $n$-dimensional Euclidean space, $\mathbb{R}=\mathbb{R}^1$, $\mathbb{R}^+$ is the set of non-negative real numbers, and $\mathbb{R}^{n\times m}$ is the real valued $n\times m$ matrix space.

We use lowercase letters to denote column vectors or scalars, uppercase letters to denote matrices. An identify matrix is $I$.  The superscript $T$ denotes the transpose of a vector or matrix. The trace of a matrix $A$ is $tr(A)$.  For symmetric matrices $A$ and $B$, $A\geq 0$ or $A>0$ denotes that $A$ is positive semi-definite or positive definite, respectively, and $A\geq B$ or $A>B$ denotes $A-B\geq 0$ or $A-B>0$, respectively. For a square matrix $A$, its determinant is $|A|$,
%its maximum (real) eigenvalue is $\lambda_{min}(A)$.
its maximum and minimum eigenvalues (when exist) are $\lambda_{max}(A)$ and $\lambda_{min}(A)$, respectively.

For a random variable $x\in \mathbb{R}^n$, its mean is $E[x]$. For a scalar function $f(x): \mathbb{R}^n\ra \mathbb{R}$, its gradient is $\nb_x f(x)=[\partial f/\partial x_1, \partial f/\partial x_2, ..., \partial f/\partial x_n]^T\in\mathbb{R}^n$, and its divergence is $div(f) = \sum_{i=1}^n\partial f/\partial x_i\in\mathbb{R}$.

%Finally, we use "$A\Rightarrow B$" as a concise form of the statement "$A$ leads to $B$".

\section{Stochastic Particle Flow Through General Homotopy}
Assume a given probability space, on which $x\in \mathbb{R}^n$ is a $n$-dimensional random variable and $z\in \mathbb{R}^d$ is a $d$-dimensional measurement of $x$. Let $p_x(x)$ denote the prior probability density function of $x$ and $p_z(z|x)$ the likelihood of a measurement $z$ conditioned on $x$. The Bayes' Theorem states that the posterior conditional density function of $x$ for a given measurement $z$, $p_x(x|z)$, is given by\footnote{To keep notations as simple as possible, in this paper we focus on one-step Bayesian estimation which can be applied to filtering or inference problems. For multi-step sequential filtering, the Bayes' Theorem is as the following \cite{DH2007,Jaz}
\[
p(x,t_k|Z_k)=p(z_k|x,t_k)p(x,t_k|Z_{k-1})/p(z_k|Z_{k-1})
\]
in which $z_k$ is the $k$-th measurement at time $t_k$, $Z_k=\{z_1,z_2,...,z_k\}$, and $p(z_k|x,t_k)$ is the probability density of measurement $z_k$ at time $t_k$ conditioned on $x$. The probability density functions $p_0(x)$, $p_1(x)$, and $h(x)$ in (\ref{pgh}) need to be replaced, respectively, with the following
\[ 
p_0(x)=p(x,t_k|Z_{k-1}), p_1(x)=p(x,t_k|Z_k), h(x)=p(z_k|x,t_k).
\]
The rest of discussion follows.
}
\bq
p_x(x|z) = \frac{p_x(x)p_z(z|x)}{p_z(z)}
\label{bayes}
\eq
in which $p_z(z)=\int_x p_z(z|x)p_x(x)dx$ is the normalization factor. Without loss of generality, it is assumed throughout this paper that all probability density functions exist, second order continuously differentiable, and are non-vanishing everywhere. Those assumptions are stronger than we need but are helpful to maintaining clarity of discussion without resorting to complex notations. %Otherwise, we restrict discussion to the supports of the density functions.
For notational simplicity, we denote
\bq p_0(x) = p_x(x), p_1(x)=p_x(x|z), h(x) = p_z(z|x).
\label{pgh}
\eq
Then $p_0(x)$ is the prior distribution and $p_1(x)$ is the posterior distribution, and $h(x)$ is the measurement likelihood. 

Let $\al(\lambda)\in \mathbb{R}, \be(\lambda)\in \mathbb{R}$ be two continuously differentiable scalar functions defined over $\lambda\in [0, 1]$ such that 
\bq
\al(0)=1, \al(1)=0, \be(0)=0, \be(1)=1.
\label{albe_bc}
\eq
For each $\lambda \in [0, 1]$, define a new conditional probability density function as the following 
\bq
p(x,\lambda) = \frac{p_0^{\al(\lambda)}(x)p_1^{\be(\lambda)}(x)}{\Gamma(\lambda)}
\label{homotopy}
\eq
in which $\Gamma(\lambda)$ is the normalization factor so that $p(x,\lambda)$ remains a probability density function. Note that 
\[    p(x,0) = p_0(x), p(x,1)=p_1(x). \]
The mapping $p(x,\lambda): \mathbb{R}^+\times [0, 1]\longrightarrow \mathbb{R}^+$ in (\ref{homotopy}) defines a homotopy from $p_0(x)$ to $p_1(x)$ for any given pair $\al(\lambda)$ and $\be(\lambda)$ satisfying (\ref{albe_bc}). By taking the natural logarithm on both sides of (\ref{homotopy}), we obtain
\bq
\log p(x,\lambda) = \al(\lambda)\log p_0(x)+ \be(\lambda) \log p_1(x) -\log \Gamma(\lambda).
\label{loghom}
\eq
In other words, $\log p(x,\lambda)$ is a $\lambda$-varying linear combination of $\log p_0(x)$ and $\log p_1(x)$.
%Recall that a major problem with particle filters is “particle degeneracy" \cite{DH2007,DJ,PS,RAG}. To mitigate this issue, particle flow methods move (change) $x$ as a function of $\lambda$, $x(\lambda)$, so that (\ref{homotopy}), or equivalently (\ref{loghom}), is always satisfied as $\lambda$ changes from $0$ to $1$. The probability density function of $x(\lambda)$ is $p(x,\lambda)$ for all $\lambda\in [0,1]$. The value of $x(\lambda)$ at $\lambda=1$ is used for estimation in problems such as filtering or Bayesian inference. 
We next derive particle flow filters based on the log-homotopy (\ref{loghom}). There exists much freedom in the choice of $\{x(\lambda), \lambda\in[0, 1]\}$ \cite{DH2015}.
For example, $x(\lambda)$ could be driven by a deterministic process as in the Exact Flow \cite{DH2015,DHN2010}, or by a stochastic process as in stochastic flows \cite{DH2013,DHN2016,DHN2018}. Generally speaking, stochastic flows have been observed numerically stable. In this paper, we focus on stochastic flows in which $x(\lambda)$ is driven by the following stochastic process
\bq
dx = f(x,\lambda)d\lambda+q(x,\lambda)dw_{\lambda}
\label{flow}
\eq
where $f(x,\lambda)\in \mathbb{R}^n$ is a drift function, $q(x,\lambda)\in \mathbb{R}^{n\times m}$ is a diffusion matrix, and $w_{\lambda}\in \mathbb{R}^m$ is a $m$-dimensional Brownian motion process in $\lambda$ with $E[dw_\lambda dw_\lambda^T ]=\sigma(\lambda)d\lambda$. The stochastic differential equation (\ref{flow}) is a standard diffusion process \cite{Jaz}. Note that $\{x(\lambda), \lambda \in [0, 1]\}$ is a stochastic process in $\lambda$, not in time. For clarity, we drop its dependence on $\lambda$ but add the dependence back when it is beneficial to emphasize its dependence on $\lambda$. Without loss of generality, we assume that $\sigma(\lambda)=I_{m\times m}$, and denote
\[ Q(x,\lambda) = q(x,\lambda)q(x,\lambda)^T \in \mathbb{R}^{n\times n}. \]
The matrix $Q(x,\lambda)=[Q_{i,j}]$ is always symmetric positive semi-definite for any $x$ and $\lambda$.

Our goal is to select $f(x,\lambda)$ and $q(x,\lambda)$ (or equivalently $Q(x,\lambda)$) such that (\ref{loghom}) is maintained  for the particle $x(\lambda)$ driven by the stochastic process (\ref{flow}) for all $\lambda\in [0, 1]$. To that end, we start with the following necessary condition.
\vs

{\sc Lemma 2.1.} \cite{DH2013} {\em  Assume that all derivatives exist and are continuous. For the particle flow $x(\lambda)$ defined in (\ref{flow}), a necessary condition for its density function to be $p(x,\lambda)$ for all $\lambda\in [0, 1]$ is that $f(x,\lambda)$ and $Q(x,\lambda)$ satisfy the following condition
\bq
\nb_x\frac{\partial \log p}{\partial\lambda} = - \nb_x div(f)-(\nb_x\nb_x^T\log p)f 
-(\nb_x^T f)(\nb_x\log p)+\nb_x[\frac{1}{2p}\sum_{i,j}\frac{\partial^2 (pQ_{i,j})}{\partial x_i\partial x_j}] 
\label{cond1}
\eq 
for all $x$ and $\lambda\in[0, 1]$. For simplicity and without causing confusion, in (\ref{cond1}) and for the rest of discussion in this paper, we omit all variables involved.
}
\vs

Since the introduction of particle flows in \cite{DH2007}, there have been steady efforts in the literature either to solve (\ref{cond1}) for a special $Q$ or to find an approximate solution for general $Q$ \cite{Daum2016}. In this paper, we focus on finding the exact solution  $f$ for arbitrary symmetric positive semi-definite matrix $Q$ as long as $Q$ is not a function of $x$. The matrix $Q$ may be a function of $\lambda$.
\vs

{\sc Theorem 2.1.} {\em Assume that
\begin{description}
\item[(A1)] $\nb_x \log p_0$ and $\nb_x \log p_1$ are linear in $x$, and
\item[(A2)] $\nb_x \nb_x^T \log p$ is non-singular for all $\lambda \in [0, 1]$.
\end{description}
Then, for any $Q$ positive semi-definite that is not a function of $x$, (\ref{cond1}) is satisfied by the following $f$
\bq
f = K_1(\nb_x\log p)+K_2(\nb_x\log h),
\label{f1}
\eq
\bq
K_1 = \frac{1}{2}Q
+\frac{\al\dot{\be}-\dot{\al} \be}{2(\al+\be)}(\nb_x\nb_x^T\log p)^{-1}(\nb_x\nb_x^T\log h) (\nb_x\nb_x^T\log p)^{-1}-\frac{\dot{\al}+\dot{\be}}{2(\al+\be)}(\nb_x\nb_x^T\log p)^{-1},
\label{K1}
\eq
\bq  
K_2 = - \frac{\al\dot{\be}-\dot{\al} \be}{\al+\be}(\nb_x\nb_x^T\log p)^{-1},
\label{K2}
\eq
where $\dot{\al} = \frac{d\al}{d\lambda}, \dot{\be}=\frac{d\be}{d\lambda}$.
}
\vs

Under the assumptions (A1) and (A2) in Theorem 2.1, (\ref{f1}) shows that the function $f$ is a linear combination of prior knowledge $\nb_x\log p_0$ and measurement likelihood information $\nb_x\log h$. The assumptions (A1) and (A2) are satisfied if $p_0$ and $p_1$ are Gaussian.
% and the measurement equation is
%\bq
%dz = Hxd\lambda + dv
%\eq
%in which $H\in R^{d\times n}$ does not depend on $x$, and $dv$ is a diffusion process in $\lambda$. 
Such Gaussian assumption as an approximation is widely adopted in many fields of studies.
\vs

The form of $f$ (\ref{f1}) is valid for any $Q$, $\al$ and $\be$ satisfying (\ref{albe_bc}) and under the assumption (A2). However, note that (\ref{cond1}), consequently (\ref{f1}), is a necessary condition only. We need to establish that $f$ is also sufficient in the sense that the density function of the flow $x$ (\ref{flow}) with $f$ defined in (\ref{f1}) indeed has the correct posterior distribution (\ref{homotopy}).
\vs

{\sc Theorem 2.2.} {\em Assume the assumptions (A1) and (A2) in Theorem 2.1, that $q$ is not a function of $x$, and that the flow $x$ is governed by the diffusion process (\ref{flow}) with $f$ defined in (\ref{f1}). Then the probability density function
of $x$ is indeed $p(x,\lambda)$ (\ref{homotopy}).
}
\vs

Theorem 2.2 establishes the correctness of the flow (\ref{f1}) in which $Q$, $\al$ and $\be$ serve as design parameters that a user can adjust for a specific problem, say, to improve filtering performance. In the remaining of this section, we examine their choices in ensuring the stability of particle flows. In the next section, we explore their choice for stiffness mitigation. 
\vs

In implementing particle flow filters, we need to numerically evaluate the stochastic differential equation (\ref{flow}) to obtain a numerical approximation of its theoretical exact solution, for given initial condition. The approximate solution in turn is used to construct estimates using the value of $x$ at $\lambda=1$. Therefore, ensuring the stability of numerically solving (\ref{flow}) is an important issue. In \cite{DD2021}, we established the stability of particle flows for the straight line log-homotopy (\ref{line_loghomo}). In this paper, we extend the analysis to the more general form of homotopy (\ref{homotopy}).

Consider the flow (\ref{flow}) with $q$ not a function of $x$. Under the assumptions (A1) and (A2) in Theorem 2.1, $f$ is linear in $x$. Therefore, we can separate out the linear and the constant terms, and re-write $f$ as 
\[ 
f = F(\lambda)x+b(\lambda)
\]
in which $F(\lambda)\df\nb_xf$ is the Jacobian matrix of $f$ and $b(\lambda)=f-F(\lambda)x$. According to (\ref{f1})-(\ref{K2}),
\[
F(\lambda) = K_1\nb_x\nb_x^T\log p+K_2\nb_x\nb_x^T\log h
\]
%\[
%=\frac{1}{2}Q(\nb_x\nb_x^T\log p)
%+\frac{\al\dot{\be}-\dot{\al} \be}{2(\al+\be)}(\nb_x\nb_x^T\log p)^{-1}(\nb_x\nb_x^T\log h) -\frac{\dot{\al}+\dot{\be}}{2(\al+\be)}I
%\]
%\[- \frac{\al\dot{\be}-\dot{\al}\be}{\al+\be}(\nb_x\nb_x^T\log p)^{-1}(\nb_x\nb_x^T\log h).
%\]
\bq
=\frac{1}{2}Q(\nb_x\nb_x^T\log p)
-\frac{\al\dot{\be}-\dot{\al} \be}{2(\al+\be)}(\nb_x\nb_x^T\log p)^{-1}(\nb_x\nb_x^T\log h) -\frac{\dot{\al}+\dot{\be}}{2(\al+\be)}I.
\label{Jacobian1}
\eq
%Since $Q$ is not a function of $x$, we may find a $q$ that is not a function of $x$ in (\ref{flow}) ($q$ could be a function of $\lambda$). 
Let $x_1(\lambda)$ and $x_2(\lambda)$ be two solutions to (\ref{flow}) with different initial conditions $x_1(0)$ and $x_2(0)$, respectively.  Then their difference
\[ 
\tilde{x}(\lambda)\df x_1(\lambda)-x_2(\lambda)
\]
satisfies
\bq
d\tilde{x} = F(\lambda)\tilde{x}d\lambda, \hs \tilde{x}(0)=\tilde{x}_0
\label{linearized}
\eq
in which the initial condition is $\tilde{x}_0\df x_1(0)-x_2(0)$.

We next examine the stability of (\ref{linearized}). If (\ref{linearized}) is stable, effects of an error in initial condition $\tilde{x}_0$ on the accuracy of $x(1)$ are likely limited. Stability of (\ref{linearized}) is therefore desirable. Note that (\ref{linearized}) is a $\lambda$-varying system. We adopt an approach of Lyapunov stability that has been proven a powerful tool for analyzing stability of time-varying systems \cite{Kha}. Following \cite{DD2021}, define $M(\lambda)=-\nb_x\nb_x^T\log p$ and consider
\bq
V(\tilde{x},\lambda)=\tilde{x}^TM(\lambda)\tilde{x}.
\label{V}
\eq  
The reason for having the negative sign in the definition of $M(\lambda)$ is that the Hessian matrix $\nb_x\nb_x^T\log p$ (assuming it exists) is typically negative definite, e.g., for unimodal distributions. With the negative sign, $M(\lambda)$ is positive definite. The next lemma describes the dynamics of $V(\tilde{x},\lambda)$ which serves as a stepping stone in the stability analysis of (\ref{linearized}).
\vs

{\sc Lemma 2.2.} {\em Assume the assumptions (A1) and (A2). Then for any given $\tilde{x}_0$, we have
\bq
dV= -\tilde{x}^T(MQM)\tilde{x}d\lambda,
\label{dV}
\eq
with initial condition 
\[
V_0\df V|_{\lambda=0} = - \tilde{x}_0^T(\nb_x\nb_x^T\log p_0 )\tilde{x}_0.
\]
}
\vs

Note that $\nb_x\nb_x^T\log p$ is the Hessian matrix of the density $p$, and thus $M$ is symmetric. The right hand side of (\ref{dV}) is always non-positive for $Q$ positive semi-definite. We know from (\ref{a_ddlogp}) in the Appendix that 
\bq
M(\lambda) =-(\al+\be)\nb_x\nb_x^T\log p_0-\be\nb_x\nb_x^T\log h.
\label{M}
\eq
The Hessian matrices $\nb_x\nb_x^T\log p_0$ and $\nb_x\nb_x^T\log h$ are typically negative definite or semi-definite, respectively. This is the case, for example, for Gaussian density functions. For the straight line homotopy, we have $\al=1-\lambda$ and $\be=\lambda$. In this case, $M(\lambda) =-\nb_x\nb_x^T\log p_0-\lambda\nb_x\nb_x^T\log h$ is always positive definite if $\nb_x\nb_x^T\log p_0$ is negative definite and $\nb_x\nb_x^T\log h$ is negative semi-definite. For general homotopy (\ref{homotopy}), we need further conditions on $\al(\lambda)$ and $\be(\lambda)$ for stability analysis.
\vs

{\sc Theorem 2.3.} {\em Assume the assumptions (A1) and (A2), and that
\begin{description}
\item[(A3)] There exists a constant symmetric positive definite matrix $M_0\in \mathbb{R}^{n\times n}$ such that $M(\lambda)\geq M_0>0, \forall \lambda \in [0, 1]$.
\end{description}
Then
\bi
\item[(1)] For any $Q\geq 0$, $\tilde{x}$ is bounded such that 
\bq
\tilde{x}^TM_0\tilde{x}\leq c, \forall \lambda \in [0, 1].
\label{bounded}
\eq
\item[(2)] If $Q>0$, $\tilde{x}$ decreases exponentially in the sense that
\bq
\tilde{x}^TM_0\tilde{x}\leq c e^{-r \lambda}, \forall \lambda \in [0, 1]
\label{exp_bound}
\eq
\ei
in which $c=\tilde{x}_0^TM(0)\tilde{x}_0\geq 0, r=\lambda_{min}(Q)\lambda_{min}(M_0)>0$.
}
\vs

It has been known that noise in a diffusion process may be exploited to stabilize a dynamic system \cite{ACW,Mao}. As an example, consider the following one-dimensional stochastic differential equation
\bq
dx(t) = axdt+q xdw(t)
\label{example}
\eq
in which $x$ is a scalar random variable, $t$ is time, $a$ and $q$ are constants, and $w(t)$ is a one-dimensional Brownian motion with $E[(dw)^2]=1$. The solution $x$ is \cite{Kha2012, Koz}
\bq
x(t) = x_0e^{(a-q^2/2)t+q w(t)}.
\label{example_sol}
\eq
Then $\lim_{t\ra\infty}x(t)= 0, a.s.$ for any $x_0$ if $q^2>2a$. Therefore, (\ref{example}) is stable if $q^2>2a$. However, (\ref{example}) is unstable if $q=0$, $a>0$. The noise term stabilizes an otherwise unstable system.

For the particle flow filters existing in the literature, the diffusion matrix $Q$ is positive semi-definite, which ensures the boundedness of $\tilde{x}$ in the sense of (\ref{bounded}). If we further select $Q$ positive definite, Theorem 2.3 (2) says that any initial error will be reduced at least exponentially. This is a global exponential stability. The form (\ref{exp_bound}) suggests that a larger $\lambda_{min}(Q)$ leads to faster rate of error reduction. However, caution needs to be exercised in the choice of $Q$. This is obvious from (\ref{example_sol}): If $q$ is too large, also noting that $w(t)$ is an unbounded Gaussian random number, the exponential term $e^{(a-q^2)t+qw(t)}$ in (\ref{example_sol}) may be sufficiently large to cause numerical instability which should be avoided.

\section{Stiffness Mitigation through General Homotopy}

When particle flow filters are implemented in practice, we typically generate a number of, say $N$, particles with initial conditions $\{x_{i}(0), i=1, 2, ..., N\}$, propagate the particles to $\{x_i(1), i=1,2, ..., N\}$ by numerically solving the flow equation (\ref{flow}), and construct an estimate(s) using $\{x_i(1),i=1,2,...,N\}$. Therefore, numerically solving (\ref{flow}) with desired accuracy is important to the implementation of particle flow filters.

It has been long known that certain deterministic or stochastic differential equations could be stiff \cite{CH,HW,Lam}, which poses challenges for their numerical evaluation. While the phenomenon of stiffness has been widely observed, a precise mathematical definition of stiffness has not been available. One approach of linking stiffness to a mathematical formulation is to use the Jacobian matrix  to describe the degree of stiffness of a  differential equation, known as the linear stability theory \cite{Lam}. Consider an abstract $n$-dimensional stochastic equation
\[
\frac{dx(t)}{dx}=f(x).
\] 
Without loss of generality, assume that $f(0)=0$. Denote
\[
F=\nb_x f |_{x=0}.
\]
Assume that all eigenvalues of $F$ have strictly negative real parts. Let $\{\lambda_i, i=1, 2,...,n\}$ be the collection of the eigenvalues of $F$ and define $\overline{\lambda},\underline{\lambda} \in \{\lambda_i, i=1, 2,...,n\}$ such that
\[
|Re\overline{\lambda}|=\max_i|Re \lambda_i|, \hs
|Re\underline{\lambda}|=\min_i|Re \lambda_i|.
\]
The stiffness ratio of $F$ is defined as
\[
R_{\textrm{stiff}} = \frac{|Re\overline{\lambda}|}{|Re\underline{\lambda}|}.
\]
The stochastic differential equation (\ref{flow}) is stiff if (1) all eigenvalues of $F$ have negative real parts and (2) its stiffness ratio $R_{\textrm{stiff}}$ is large. The stiffness ratio $R_{\textrm{stiff}}$  describes the degree of stiffness. It needs to be pointed out this is a heuristic measure of stiffness: A large $R_{\textrm{stiff}}$ may (often does) lead to a stiff differential equation, but this is not guaranteed. A linear system can be constructed  with a large $R_{\textrm{stiff}}$ (e.g., due to scaling) but is not stiff \cite{Lam}. 

Such a definition links stiffness to the condition number of the matrix $F$: A large $R_{\textrm{stiff}}$ indicates the condition number of $F$ is bad (or large) and vice versa. We next examine the stiffness of the stochastic equation (\ref{flow}) by considering the condition number of its Jacobian matrix $F$, and how to choose $\al(\lambda)$ and $\be(\lambda)$ to mitigate its stiffness. Direct application of the stiffness ratio may not be meaningful: The Hessian matrix $\nb_x\nb_x^T\log h$ is singular if the dimension of measurement is lower than that of the state, i.e., $d<n$, which is often the case in practice. The diffusion matrix $Q$ could be singular also \cite{DHN2016,DHN2018}. In these cases, $F$ may have an eigenvalue zero and thus its condition number would be infinite, regardless the choice of $\al, \be$. However, a closer look of the form of $F$ in (\ref{Jacobian1}) reveals that the condition number of $\nb_x\nb_x^T\log p$ plays a critical role in affecting the condition number of $F$ (assuming $\al, \be$ are well behaved over $[0, 1]$). The inversion of a Hessian matrix $\nb_x\nb_x^T\log p$ is expected to lead to a large stiffness ratio if its condition number is large. Generally speaking, the condition number of a matrix should be small to ensure numerical stability \cite{Gen}. Henceforth, we next seek to choose $\al,\be$ so that the condition number of $\nb_x\nb_x^T\log p$ is as small as possible. The condition number of a matrix is always greater than or equal to $1$ \cite{Gen}. Therefore, it is feasible to minimize the condition number of $\nb_x\nb_x^T\log p$ if the parameter space is properly defined. 

For the matrix $M(\lambda)$ defined in (\ref{M}), its condition number is defined as
\[
\kappa_{\nu}(M)=\left\{
\ba{ll} ||M||_\nu||M^{-1}||_\nu, & \textrm{if $M$ is nonsingular,} \\
\infty, & \textrm{if $M$ is singular},
\ea\right.
\]
in which $||\cdot||_\nu$ is a matrix norm. The condition number of $M$ is only affected by $\be/(\al+\be)$. Without loss of generality, we normalize $\al$, $\be$ by setting $\al+\be=1$. In this case, $\dot{\al}+\dot{\be}=0$ and the Jacobian matrix (\ref{Jacobian1}) becomes
\bq
F=\frac{1}{2}Q(\nb_x\nb_x^T\log p)
-\frac{\dot{\be}}{2}(\nb_x\nb_x^T\log p)^{-1}(\nb_x\nb_x^T\log h)
\label{Jacobian2}
\eq
and 
\[
M=-\nb_x\nb_x^T\log p=-\nb_x\nb_x^T\log p_0-\be\nb_x\nb_x^T\log h.
\]
Our goal is to choose $\be$ to minimize the condition number of the matrix $M$ (or equivalently $\nb_x\nb_x^T\log p$). 

%One way of approaching this problem is through Euler-Lagrange equation. 
Traditional approaches to improving condition numbers typically use matrix analysis techniques such as scaling \cite{BM}. By exploiting the special structure of $M$, we formulate the problem of choosing $\be(\lambda)$ to minimize the condition number of $M$ as an optimal control problem over $[0, 1]$, described as the following:
\bq
\frac{d\be}{d\lambda}=u(\lambda), 
\label{state}
\eq
\bq
\be(0)=0, \be(1)=1,
\label{BV}
\eq
\bq
J(\be,u)=\int_0^1[\frac{1}{2}u^2+\mu\kappa_{\nu}(M)]d\lambda,
\label{obj}
\eq
in which $\mu\geq 0$ is a user-defined scalar weight. Our goal is to choose $u(\lambda)\in\mathbb{R}$ such that the objective function $J$ is minimized. In (\ref{obj}), the first term inside the integrand is interpreted as energy consumed. Therefore, this problem formulation represents a trade-off between energy needed and small condition number.
\vs

{\sc Theorem 3.1.} {\em For the optimal control problem (\ref{state})-(\ref{obj}), the optimal solution $\be^*$ of $\be$ is the solution of
\bq
\frac{d^2\be^*}{d\lambda^2}=\mu\frac{\partial \kappa_{\nu}(M)}{\partial \be}|_{\be=\be^*},
\label{opti_b}
\eq
\bq
\be^*(0)=0, \be^*(1)=1.
\label{opti_bv}
\eq
}

The optimal solution (\ref{opti_b})-(\ref{opti_bv})
is a second order two-point boundary value problem with Dirichlet boundary condition (\ref{opti_bv}) \cite{SB}.  The initial condition $\be^*(0)=0$ is given as a constraint in (\ref{albe_bc}). We need to choose  $\dot{\be}^*(0)\df\frac{d\be^*}{d\lambda}|_{\lambda=0}$ such that $\be^*(1)=1$.
In (\ref{opti_b}), the partial differential should be understood as subdifferential if the matrix norm is not continuously differentiable, e.g., $L_1$ norm. Two-point boundary value problems often arise in optimal control, astrodynamics, physics, as well as many other applications \cite{BH}. The study of two-point boundary value problems is well established. There are a number of methods available to solve two-point boundary value problems such as the theory of low and upper solutions, the (simple or multiple) shooting methods, the difference methods, and the variational methods \cite{DCH,SB}.
\vs

{\sc Remark 3.1.} Consider a special case of $\mu=0$. This is a problem of minimum energy control. The optimal solution equation (\ref{opti_b}) becomes
\[
\frac{d\be^*}{d\lambda}=0, \be^*(0)=0, \be^*(1)=1.
\]
Its solution is the straight line $\be^*=\lambda$ and consequently $\al^*=1-\lambda$, which recovers the linear convex combination that has been used to construct homotopy in particle flow filters in the literature \cite{DH2015,DHN2018}. Note that $\dot{\be}=u$. This case corresponds to minimizing the contribution of the second term to the Jacobian matrix $F$ in (\ref{Jacobian2}). If $\nb_x\nb_x^T\log p$ is ill-conditioned (i.e., with a large condition number), the accuracy of its inverse is expected to be less. Minimizing $u$ would minimize its contribution to the Jacobian matrix $F$. 
\vs

{\sc Remark 3.2.} The Hessian matrix $\nb_x\nb_x^T\log p$ is typically negative definite for practical problems, e.g., for unimodal Gaussian density functions. In this case, $M$ is positive definite. The nuclear norm of a positive (semi-)definite matrix is equal to its trace: $||M||_*=tr(M)$. Therefore, $\kappa_*(M) = tr(M)tr(M^{-1})$. Note that $M$ is linear in $\be$ and $\partial M/\partial \be = - \nb_x\nb_x^T\log h$. Consequently, noting the equality for differentiating an inverse of a matrix in Lemma A.1, (\ref{opti_b}) becomes
\[
\frac{d^2\be^*}{d\lambda^2}=-\mu [tr(\nb_x\nb_x^T\log h)tr(M^{-1})+ tr(M)tr(M^{-1}\nb_x\nb_x^T\log h M^{-1})]
\]
\bq
=-\mu [tr(\nb_x\nb_x^T\log h)tr(M^{-1})+ tr(M)tr(M^{-2}\nb_x\nb_x^T\log h)]
\label{kappa_L1}
\eq
in which we used the property $tr(AB)=tr(BA)$ in the last step.
%in which $N=\nb_x\nb_x^T\log h$.
\vs

{\sc Remark 3.3.} For the $L_2$ norm, or the spectral norm, $||M||_2=\lambda_{\max}(M)$ and $||M^{-1}||_2=\lambda_{\min}^{-1}(M)$ for a positive definite $M$. Assume that the largest and the smallest eigenvalues of $M$ are unique, and $v_{\max}$ and $v_{\min}$ are corresponding unit eigenvectors, respectively. In this case, \cite{PP}
\[
\frac{\partial ||M||_2}{\partial \be}=v_{\max}^T(\frac{\partial M}{\partial\be}) v_{\max} = -v_{\max}^T(\nb_x\nb_x^T\log h) v_{\max},
\]
\[
\frac{\partial ||M^{-1}||_2}{\partial \be}=-\lambda_{\min}^{-2}(M)v_{\min}^T(\frac{\partial M}{\partial\be}) v_{\min} = \lambda_{\min}^{-2}(M)v_{\min}^T(\nb_x\nb_x^T\log h) v_{\min}.
\]
Therefore,
\[
\frac{\partial \kappa_2(M)}{\partial \be}=
\frac{\partial ||M||_2}{\partial \be}||M^{-1}||_2+||M||_2\frac{\partial ||M^{-1}||_2}{\partial \be}
\]
\[
=-\lambda_{\min}^{-1}(M)v_{\max}^T(\nb_x\nb_x^T\log h) v_{\max}+\lambda_{\max}(M)\lambda_{\min}^{-2}(M)v_{\min}^T(\nb_x\nb_x^T\log h) v_{\min}.
\]
\vs

Of course, minimizing the condition number of $M=-\nb_x\nb_x^T\log p$ is different from improving that of the Jacobian matrix $F$. However, improving the condition number of $F$ could be incorporated in the optimal control problem formulation as a constraint on $\be$. We can easily modify the optimal solution to ensure no degradation over the straight line log-homotopy. For notational clarity, we re-write $F$ in (\ref{Jacobian2}) as $F(\lambda,\be)$ to highlight its dependence on $\be$. The straight line $\be_l=\lambda$ is the baseline. We modify the optimal $\be^*(\lambda)$ in the following way, assuming the modification is only needed at a finite number of points of $\lambda$,
\[
\be^*_{mod}=\left\{
\ba{ll} 
\be^*(\lambda), & \textrm{if } \kappa_{\nu}(F(\lambda,\be^*))\leq \kappa_{\nu}(F(\lambda,\be_l))), \\
\be_l, & \textrm{if } \kappa_{\nu}(F(\lambda,\be^*))> \kappa_{\nu}(F(\lambda,\be_l)).
\ea
\right.
\]
Then we have $\kappa_{\nu}(F(\lambda,\be^*_{mod}))\leq \kappa_{\nu}(F(\lambda,\be_l))$ for all $\lambda \in [0, 1]$. In other words, the new $\be^*_{mod}$ leads to an $F$ with condition number smaller than or equal to that for the baseline log-homotopy for all $\lambda \in [0, 1]$.
\vs

A matrix norm $||\cdot||_\nu$ is said monotone if $||B||_\nu \leq ||A||_\nu$ for any positive definite matrices $A$ and $B$ satisfying $0\leq B\leq A$. The $L_2$ and the Frobenius norms are monotone, and the nuclear norm is not. For a monotone norm, the following Theorem 3.2 states that the optimal $\be^*$ is always non-negative if the condition number of $-\nb_x\nb_x^T\log h$ is no greater than that of $-\nb_x\nb_x^T\log p_0$.
\vs

{\sc Theorem 3.2.} {\em Consider a monotone norm $||\cdot||_\nu$. Assume that the optimal solution $\be^*$ is continuously differentiable, and that both $-\nb_x\nb_x^T\log p_0$ and $-\nb_x\nb_x^T\log h$ are positive definite. If  $\kappa_\nu(\nb_x\nb_x^T\log h)\leq \kappa_\nu(\nb_x\nb_x^T\log p_0)$, we must have
$\be^*(\lambda)\geq 0$ for all $\lambda \in [0, 1]$.
}
\vs

The optimal solution $\be^*$ is obtained by minimizing the condition number of $M=(-\nb_x\nb_x^T\log p_0)+\be(-\nb_x\nb_x^T\log h)$. If the condition number of $-\nb_x\nb_x^T\log h$ is less than that of $-\nb_x\nb_x^T\log p_0$, then $\be^*$ should place more weight on $-\nb_x\nb_x^T\log h$ to reduce the overall condition number of $M$, which can only be achieved with a non-negative $\be^*$ according to Lemma A.3 (note that both $-\nb_x\nb_x^T\log p_0$ and $-\nb_x\nb_x^T\log h$ are positive definite): A negative $\be^*$ would increase the condition number of $M$ according to (\ref{a_cn_leq}). It should be pointed out that Theorem 3.2 provides a sufficient condition only. Further research is needed to fully understand the properties of $\be^*$.
\vs

Flow stability is important to maintaining numerical stability in the implementation of particle flow filters. We need to verify that $\be^*$ leads to a stable flow as discussed in Theorem 2.3. As a sufficient condition, we may check that the assumption (A3) is satisfied for the optimal $\be^*$. If we know that $\be^*(\lambda)\geq 0$ for all $\lambda \in [0, 1]$, $-\nb_x\nb_x^T\log p_0>0$, and $-\nb_x\nb_x^T\log h\geq 0$, the assumption (A3) is satisfied with $M_0=-\nb_x\nb_x^T\log p_0$. Theorem 3.2, together with Theorem 2.3, provides a sufficient condition for the stability of particle flows without explicitly checking the assumption (A3). 

Comparing  with using the straight line log-homotopy, the properties of the particle flow filters derived in this paper is remarkably consistent with those in \cite{DD2021} in terms of particle distributions and flow stability, which is not surprising because $\be=\be^*$ represents a diffeomophism of the straight line $\be=\lambda$. Solving the two-point boundary problem (\ref{opti_b}) is a small price we pay for stiffness mitigation. More importantly, the results point out that it could be beneficial to consider general forms of homotopy, which warrants further research. One possibility is to explicitly include a term of stiffness mitigation in the derivation of particle flows.

\section{A Numerical Example}

We use a numerical example to illustrate the effectiveness of the proposed approach to stiffness mitigation. The scenario is a modification from a 2D experiment in \cite{MDD}, and is graphically described in Figure 1: There are two passive infrared sensors located at (3.5, 0) and (-3.5, 0). The truth location of a stationary target is (4,4).

\begin{figure}[h]
	\centering
	\includegraphics[width=4in]{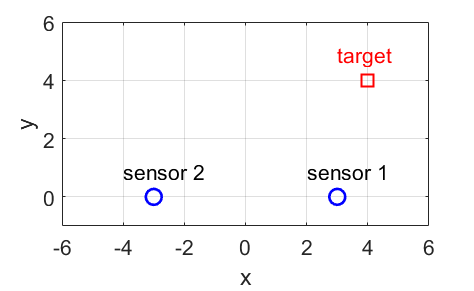}
	\caption{Experiment scenario}
	\label{fig:scenaro}
\end{figure}

For sensor $i, i=1, 2$, its measurement is the angle of the target relative to the sensor as the following
\[
h_i  = \arctan(\frac{y_t-y_i}{x_t-x_i}) + v_i, i = 1, 2,
\]
in which $(x_t,y_t)$ is the coordinate of the target, $(y_i,x_i)$ is the coordinate of sensor $i$, and $v_i$ is the (Gaussian) measurement noise. The overall measurement equation is
\[
h = \left[ \ba{c} h_1 \\ h_2 \ea \right].
\]
The prior distribution is Gaussian with mean 
\[
x_{prior} = \left[ \ba{c} 3.0 \\ 5.0 \ea \right]
\]
and covariance matrix
\[
P_{prior} = \left[ \ba{cc} 1000.0 & 0 \\ 0 & 2.0 \ea \right].
\]
The prior knowledge about the location of the target is quite off. For the measurement, the measurement noise $[v_1, v_2]^T$ is Gaussian with zero mean and covariance matrix
\[
R = \left[ \ba{cc} 0.04 & 0 \\ 0 & 0.04 \ea \right].
\]
We take a measurement as a sample automatically generated in Matlab
\[
z=\left[ \ba{c} 0.4754 \\ 1.1868 \ea \right].
\]

Other parameters are $\mu=0.2$ and 
\[
Q=\left[ \ba{cc} 4.0 & 0 \\ 0 & 0.4 \ea \right].
\] 
Experiments were carried out using Matlab R2018b installed on a regular laptop computer. Matlab's built-in solver ode45 is called to solve (\ref{opti_b}) for given initial conditions.
The nuclear norm based condition number (\ref{kappa_L1}) is used. The simple bisection method, a special case of the shooting method, is used to find $\dot{\be}^*(0)$ such that $\be^*(1)=1$. %The initial interval for $\dot{\be}^*(0)$ is $[-10, -20]$. 
The objective function $J$ for the straight line log-homotopy is 4.0, while the objective function corresponding to the optimal $\be^*(\lambda)$ is 3.4.

Figure 2 shows comparison of the optimal $\be^*(\lambda)$ with the straight line log-homotopy $\be(\lambda)=\lambda$. Figure 2 (d) shows reduction of the stiffness of the Jacobian matrix $F$ for the flow equation (\ref{flow}). Note that $R_{\textrm{stiff}}$ is plotted in log-scale to aid visualization. The reduction in the stiffness ratio is modest for this example. However, the impact on reducing estimation errors is significant as shown in Table~\ref{tab:table_homo}.

\begin{figure}[h]
	\centering
	$\ba{cc}
		\includegraphics[width=3in]{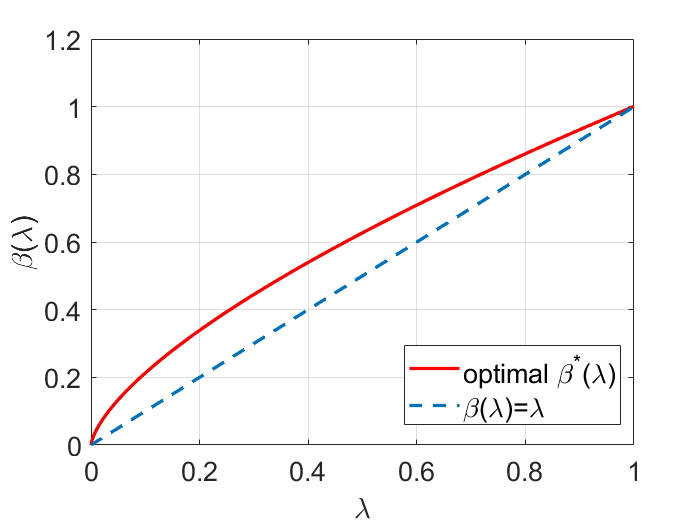}&
		\includegraphics[width=3in]{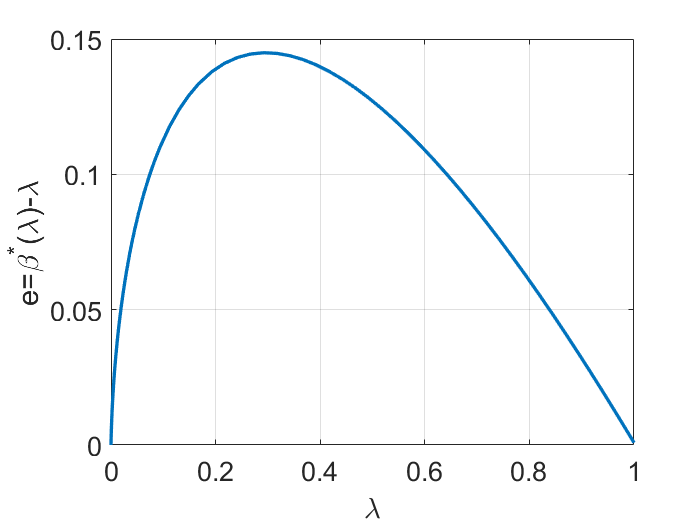} \\
		(a) \textrm{ } \be(\lambda) & (b) \textrm{ } e=\be^*(\lambda)-\lambda \\
		\includegraphics[width=3in]{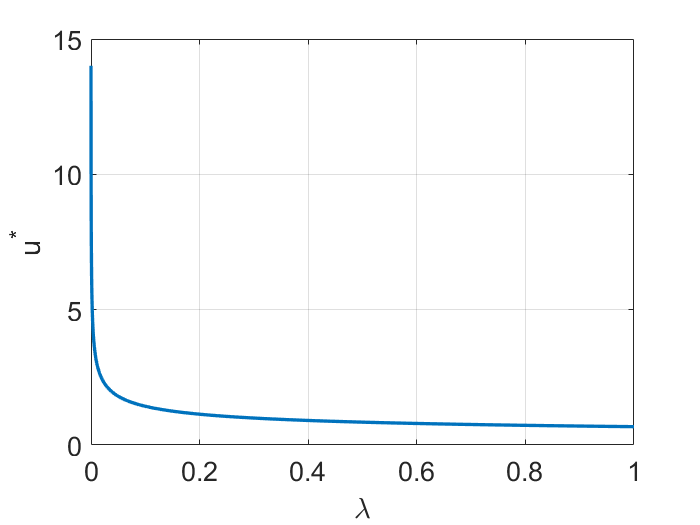}&
		\includegraphics[width=3in]{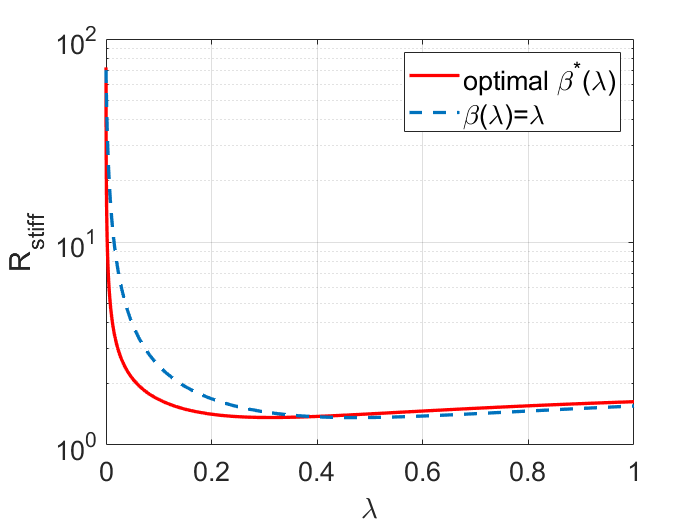} \\
		(c)  \textrm{ } u^*(\lambda) & (d) \textrm{ } R_{\textrm{stiff}}
	\ea$
	\caption{Comparison of the optimal $\be^*(\lambda)$ with the straight line $\be(\lambda)=\lambda$.}
	\label{fig:comparison}
\end{figure}

We performed 20 Monte Carlo (MC) runs to evaluate the filtering accuracy for both the baseline homotopy and the optimal homotopy. For each run, the number of particles is 50. Same sequence of Brownian motion samples are used in the numerical evaluation of the flow (\ref{flow}) for both the straight line and the optimal $\be^*(\lambda)$, which has been known as the scheme of Common Random Numbers (CRN) in Monte Carlo methods \cite{BFS}. CRN allows that the evaluation be focused on the comparison of numerical accuracy by taking out the underlying randomness. The estimation accuracy is measured by the Mean Square Error (MSE) of the estimate and the trace of error covariance matrix $tr(P)$. The results are shown in Table~\ref{tab:table_homo} for each MC run and overall average. The optimal homotopy indeed improves filtering accuracy in terms of reduced MSE and covariance matrix. In Table~\ref{tab:table_homo}, the subscription $\be$ refers to the case of the straight line log-homotopy $\be(\lambda)=\lambda$.

\begin{table}[h!]
	\renewcommand{\arraystretch}{1.3}
	\begin{center}
		\caption{Performance Comparison}
		\vs
		\label{tab:table_homo}
		\begin{tabular}{|c|c|c|c|c|}
			\hline
			 MC index  & $MSE_{\be}$ & $MSE_{\be^*}$ & $tr(P_\be)$ & $tr(P_{\be^*})$ \\ \hline
1 & 11.658 & 9.0751 & 1172.7 & 867.37  \\ \hline
2 & 20.287 & 14.558 & 2099.7 & 1284.6  \\ \hline
3 & 9.7335 & 6.3453 & 1106.4 & 644.72  \\ \hline
4 & 12.119 & 9.7331 & 1520.5 & 1053.2  \\ \hline
5 & 9.7364 & 6.4886 & 814.15 & 604.49  \\ \hline
6 & 8.2546 & 4.8693 & 1625.4 & 1045.7  \\ \hline
7 & 14.01  & 10.83  & 1689.8 & 1195.2  \\ \hline
8 & 17.743 & 14.085 & 1804   & 1245.8  \\ \hline
9 & 17.443 & 12.221 & 2008.8 & 1306.8  \\ \hline
10 & 12.403 & 8.9431 & 1106.5 & 720.68 \\ \hline
11 & 19.18  & 13.669 & 2693.5 & 1749.6 \\ \hline
12 & 15.668 & 12.304 & 1792.9 & 1311.6 \\ \hline
13 & 19.151 & 13.839 & 1926   & 1280.7 \\ \hline
14 & 15.285 & 9.8745 & 1985.2 & 1132.9 \\ \hline
15 & 9.1811 & 5.4055 & 862.67 & 537.48 \\ \hline
16 & 6.0954 & 2.9439 & 919.84 & 626.66 \\ \hline
17 & 9.7702 & 6.4905 & 1217.2 & 822.45 \\ \hline
18 & 8.3131 & 5.475  & 1101.8 & 801.86 \\ \hline
19 & 8.4433 & 6.7426 & 1080.3 & 840.67 \\ \hline
20 & 20.456 & 15.614 & 2176   & 1503.8 \\ \hline
average & 13.246 & 9.4754 & 1535.2 & 1028.8 \\ \hline
		\end{tabular}
	\end{center}
\end{table}

\section{Conclusions}
The straight line log-homotopy has been almost exclusively used in  the derivation of particle flow filters in the literature. While the straight line log-homotopy has several attractive features such as easy to implement and leads to robust solutions, one natural question is: Is it beneficial to consider other forms of homotopy? In this paper, we seek answers to this question by considering a general class of $\lambda$-varying linear log-homotopy (\ref{homotopy}). The form of homotopy (\ref{homotopy}) by no means represents all possible forms of homotopy although it is certainly more general than the straight line log-homotopy. For (\ref{homotopy}), we derived particle flow filters, established the correctness of the flow by showing that the particle flows indeed have desired posterior probability distributions, and derived conditions under which the flows are stable. In Section 3, we formulated the problem of finding an effective homotopy as an optimal control problem. This novel approach explicitly exploits the special structure of the stochastic differential equations in particle flow filters. By applying the Pontryagin maximum principle, we showed that the optimal solution can be obtained by solving a second order two-point boundary value problem with Dirichlet condition. This is a one-dimensional two-point boundary value problem in a rather simple form. The solution can be obtained efficiently using standard algorithms such as the shooting method. The results on the choices of homotopy provides partial answers to relevant open problems in \cite{DH2016}. In particular, the results demonstrated that it could be beneficial to consider a more general form of homotopy, which potentially opens a new direction for research, e.g., stiffness mitigation through general homotopy or transport beyond the linear form considered in this paper. The numerical example in Section 4 demonstrated that the optimal homotopy indeed reduces the stiffness ratio of the stochastic flow equation (\ref{flow}), which in turn reduces the mean square error of the particle flow filters.

The results presented in this paper represent the first step in constructing particle flows to further improve particle flow filters and gain theoretical understanding. One immediate step is to examine other forms of homotopy for stiffness mitigation, including those proposed in \cite{DH2016} or in optimal transport theory. A comprehensive list of open problems for particle flow filters is available in \cite{DH2016}.
%Our formulation of choosing homotopy to optimize the condition number of the Hessian matrix of the posterior distribution is novel. In fact, investigation of the Hessian matrix has been a subject of an active research in optimal transport. By exploiting freedom in construction, we can improve the robustness of optimal transport solutions, which may potentially leads to new techniques for machine learning.
\vs

%\section{Acknowledgment}
\newpage

\begin{center}
	{\Large \bf Appendix}\\[2mm]
\end{center}

We first state several lemmas that will be needed in the proofs.
\vs

\noindent
{\sc \bf Lemma A.1.} \cite{Gen, PP} {\em Let $A(\theta)\in \mathbb{R}^{n\times n}$ be a parameterized non-singular square matrix and differentiable with respect to the parameter $\theta\in \mathbb{R}$ in a neighborhood of a given point. Then
\[   \frac{dA^{-1}(\theta)}{d\theta} = - A^{-1}(\theta)(\frac{dA(\theta)}{d\theta})A^{-1}(\theta).
\]
}
\vs

\noindent
{\sc \bf Lemma A.2.} \cite{Bel,Gron} (Gronwall-Bellman Inequality) {\em 
Let $u(t)\in \mathbb{R}$ be a differentiable, positive function and $\al(t)\in \mathbb{R}$ a continuous function. If the derivative of $u(t)$ satisfies
\[
\frac{d u(t)}{dt}\leq \al(t)u(t), \forall t\geq t_0
\]
then
\[  u(t)\leq u(t_0)e^{\int_{t_0}^t\al(s)ds}, \forall t \geq t_0.
\]}
In particular, if $\al(t)=-\al$, $\al>0$ a constant,
\[u(t)\leq u(t_0)e^{-\al(t-t_0)}, \forall t \geq t_0.
\]
\vs

\noindent
{\sc \bf Lemma A.3.} {\em Let $A$ and $B$ be two symmetric positive definite matrices of a given dimension, $d_1\geq 0$ and $d_2\geq 0$ are scalars satisfying $A-d_2B >0$. Consider a monotone norm $||\cdot||_\nu$, i.e., $||B||_\nu\leq ||A||_\nu$ if $0\leq B\leq A$. If $\kappa_{\nu}(B)\leq \kappa_{\nu}(A)$, we have
\bq
\kappa_{\nu}(A+d_1 B) \leq \kappa_{\nu}(A-d_2B).
\label{a_cn_leq}
\eq

Furthermore, if the strict inequality holds $\kappa_{\nu}(B) < \kappa_{\nu}(A)$ and $d_1+d_2>0$, we have
\bq
\kappa_{\nu}(A+d_1B)< \kappa_{\nu}(A-d_2 B).
\label{a_cn_less}
\eq
}

Note that the inequality (\ref{a_cn_leq}) may no longer hold if $\kappa_{\nu}(B)>\kappa_{\nu}(A)$. A counter example would be  $d_1=1$, $d_2=0$, $A=I$, $\kappa_2(A)=1$, and any $B$ such that $\kappa_2(A+B)>1$. 

{\bf Proof}: For a monotone norm $||\cdot||_\nu$, it was shown in \cite{MO1969,MO1973} that
\bq
\kappa_{\nu}(A+B) \leq \max\{\kappa_{\nu}(A), \kappa_{\nu}(B)\}
\label{a_cn_1}
\eq
for any positive definite matrices $A$ and $B$. If $d_1\geq  0$ and $\kappa_{\nu}(B)\leq \kappa_{\nu}(A)$, it directly follows from  (\ref{a_cn_1}) that
\bq
\kappa_{\nu}(A+d_1B) \leq \max\{\kappa_{\nu}(A), \kappa_{\nu}(d_1B)\} = \kappa_{\nu}(A).
\label{a_cn_2}
\eq

We next prove $\kappa_{\nu}(A)\leq \kappa_{\nu}(A-d_2B)$ for $d_2>0$. Let $C=A-d_2 B>0$. Then
$A=C+d_2B$. Define % $U_B=||B||_\nu^{-1}B$, $U_C=||C||_\nu^{-1}C$, and $\theta=||C||_\nu[||C||_\nu+d_2||B||_\nu]^{-1}$. 
\[
U_B=\frac{B}{||B||_\nu}, U_C=\frac{C}{||C||_\nu},  \theta=\frac{||C||_\nu}{||C||_\nu+d_2||B||_\nu}. 
\]
We have $||U_B||_\nu=1, ||U_C||_\nu=1$, and $0<\theta<1$. By the convexity of matrix inversion for positive definite matrices, \cite{MO1969}
\[
[\theta U_C+(1-\theta)U_B]^{-1}\leq \theta U_C^{-1}+(1-\theta)U_B^{-1}
\]
which gives, also using properties of a matrix norm,
\bq
||[\theta U_C+(1-\theta)U_B]^{-1}||_\nu\leq \theta ||U_C^{-1}||_\nu+(1-\theta)||U_B^{-1}||_\nu.
\label{a_cn_3}
\eq
On the other hand,
%\[
%\kappa_{\nu}(B) =||U_B||_\nu ||U_B^{-1}||_\nu = || U_B^{-1}||_\nu
%\]
%and 
\[
\kappa_{\nu}(A)=\kappa_{\nu}(C+d_2B)=||\theta U_C+(1-\theta)U_B||_\nu ||[\theta U_C+(1-\theta)U_B]^{-1}||_\nu
\]
\[
\leq [\theta ||U_C||_\nu+(1-\theta)||U_B||_\nu]||[\theta U_C+(1-\theta)U_B]^{-1}||_\nu
\]
\[ = ||[\theta U_C+(1-\theta)U_B]^{-1}||_\nu
\]
\[
\leq \theta ||U_C^{-1}||_\nu+(1-\theta)||U_B^{-1}||_\nu
\]
according (\ref{a_cn_3}). Therefore,
the assumption $\kappa_{\nu}(B)\leq\kappa_{\nu}(A)$ gives
\[
||U_B^{-1}||_\nu =\kappa_{\nu}(B)\leq\kappa_{\nu}(A) \leq \theta ||U_C^{-1}||_\nu+(1-\theta)||U_B^{-1}||_\nu.
\]
After consolidating terms, the previous inequality leads to
\[
||U_B^{-1}||_\nu\leq ||U_C^{-1}||_\nu
\]
or equivalently
\[
\kappa_{\nu}(B)=||U_B^{-1}||_\nu \leq ||U_C^{-1}||_\nu = \kappa_{\nu}(C).
\]
Using the previous relationship and
by applying (\ref{a_cn_1}) to $\kappa_{\nu}(A)$, we get
\[
\kappa_{\nu}(A)=\kappa_{\nu}(C+d_2 B)\leq \max\{\kappa_{\nu}(C), \kappa_{\nu}(B)\}=\kappa_{\nu}(C) =\kappa_{\nu}(A-d_2 B),
\]
which obviously also holds if $d_2=0$. The combination of the previous inequality with (\ref{a_cn_2}) gives (\ref{a_cn_leq}).

The proving process shows that (\ref{a_cn_less}) holds if we have strict inequality $\kappa_{\nu}(B) < \kappa_{\nu}(A)$ and at least one of $d_1$ and $d_2$ is not zero, i.e., $d_1+d_2>0$.
Q.E.D.
\vs

\noindent
{\sc \bf Proof of Theorem 2.1} 

First, if $Q$ is not a function of $x$ (could be a function of $\lambda$), we have
\[
\frac{\partial (\log p Q_{i,j})}{\partial x_j}
=\frac{1}{p}\frac{\partial (pQ_{i,j})}{\partial x_j}
\]
and 
\[
\frac{\partial^2 (\log p Q_{i,j})}{\partial x_i \partial x_j}
=\frac{1}{p}\frac{\partial^2 (pQ_{i,j})}{\partial x_i\partial x_j}-\frac{1}{p^2}\frac{\partial p}{\partial x_i}\frac{\partial p}{\partial x_j}Q_{i,j}
=\frac{1}{p}\frac{\partial^2 (pQ_{i,j})}{\partial x_i\partial x_j}-\frac{\partial \log p}{\partial x_i}\frac{\partial \log p}{\partial x_j}Q_{i,j}
\]
or equivalently after rearranging terms
\[
\frac{1}{p}\frac{\partial^2 (pQ_{i,j})}{\partial x_i\partial x_j}
=\frac{\partial^2 (\log p Q_{i,j})}{\partial x_i \partial x_j}+
\frac{\partial \log p}{\partial x_i}\frac{\partial \log p}{\partial x_j}Q_{i,j}.
\]
Consequently,
\[
\sum_{i,j}\frac{1}{p}\frac{\partial^2 (pQ_{i,j})}{\partial x_i\partial x_j}
=\sum_{i,j}\frac{\partial^2 (\log p Q_{i,j})}{\partial x_i \partial x_j}+\sum_{i,j}
\frac{\partial \log p}{\partial x_i}\frac{\partial \log p}{\partial x_j}Q_{i,j}
\]
\bq
=tr(\nb_x\nb_x^T\log pQ)+(\nb_x^T\log p)Q(\nb_x\log p).
\label{a_temp1}
\eq
The last equation uses the fact that $Q$ is symmetric, $Q_{i,j}=Q_{j,i}$. 
Under the assumption (A1), $\nb_x\nb_x^T\log p$ is not a function of $x$. The first term of (\ref{a_temp1}) is not a function of $x$. By taking the gradient with respect to $x$ on both sides of (\ref{a_temp1}), we obtain \cite{CL}
\bq
\nb_x[\frac{1}{2p}\sum_{i,j=1}^n\frac{\partial^2 (pQ_{i,j})}{\partial x_i\partial x_j}] =(\nb_x\nb_x^T\log p)Q(\nb_x\log p).
\label{a_eqQ}
\eq
%The equality (\ref{a_eqQ}) was derived in \cite{CL} and is re-derived here for self-containedness.
Second, note that $p_1=p_0 h/\Gamma(1)$ according to (\ref{bayes}). Therefore,
\bq
\log p_1 = \log p_0+\log h - \log \Gamma(1).
\label{a_logp1}
\eq
Substituting this form of $\log p_1$ into (\ref{loghom}),  we know that
\bq
\log p = (\al+\be)\log p_0+\be \log h -\log \Gamma(\lambda)-\log \Gamma(1).
\label{a_loghom}
\eq
Taking the gradient with respect to $x$ on both sides of (\ref{a_loghom}) gives
\bq
\nb_x\log p = (\al+\be)\nb_x\log p_0 + \be\nb_x\log h.
\label{a_dlogp_1}
\eq
Solving for $\nb_x\log p_0$ from the previous equation, we get
\bq
\nb_x\log p_0 = \frac{1}{\al+\be}\nb_x\log p-\frac{\be}{\al+\be}\nb_x\log h.
\label{a_dlogp0}
\eq
Taking the gradient again on both sides of (\ref{a_dlogp_1}), we obtain
\bq
\nb_x\nb_x^T\log p=(\al+\be)\nb_x\nb_x^T\log p_0+\be\nb_x\nb_x^T\log h.
\label{a_ddlogp}
\eq
Furthermore, by taking the derivative with respect to $\lambda$ on both sides of (\ref{a_loghom}), we have
\[
\frac{\partial \log p}{\partial \lambda}
%= \dot{\al}\log p_0+\dot{\be}\log p_1-\frac{d\log \Gamma(\lambda)}{d\lambda}
=(\dot{\al}+\dot{\be})\log p_0+\dot{\be}\log h-\frac{d\log \Gamma(\lambda)}{d\lambda}.
\]
We further take the gradient with respect to $x$ on both sides of the previous equation to eliminate the last term of normalization factor and obtain 
\bq
\nb_x[\frac{\partial \log p}{\partial \lambda}]= (\dot{\al}+\dot{\be})\nb_x\log p_0 +\dot{\be}\nb_x\log h.
\label{a_d2logpdlam_1}
\eq
Combining the form of $\nb_x\log p_0$ (\ref{a_dlogp0}) with the previous equation (\ref{a_d2logpdlam_1}), we arrive at 
\bq
\nb_x[\frac{\partial \log p}{\partial \lambda}]
=\frac{\dot{\al}+\dot{\be}}{\al+\be}\nb_x\log p +\frac{\al\dot{\be}-\dot{\al} \be}{\al+\be}\nb_x\log h.
\label{a_d2logpdlam}
\eq
Substituting both (\ref{a_eqQ}) and (\ref{a_d2logpdlam}) into the condition (\ref{cond1}) gives
\[
\frac{\dot{\al}+\dot{\be}}{\al+\be}\nb_x\log p +\frac{\al\dot{\be}-\dot{\al}\be}{\al+\be}\nb_x\log h = - \nb_x div(f)-(\nb_x\nb_x^T\log p)f 
-(\nb_x^T f)(\nb_x\log p)
\]
\bq
+(\nb_x\nb_x^T\log p)Q(\nb_x\log p).
\label{a_cond2}
\eq

For given $Q$ positive semi-definite that is not a function of $x$, we next find $f$ that satisfies (\ref{a_cond2}). Toward that end,
we consider $f$ as a linear combination of $\nb_x\log p$ and $\nb_x\log h$:
\bq
f = K_1\nb_x\log p +K_2\nb_x\log h.
\label{a_f}
\eq
Under the assumption (A1), $f$ is linear in $x$, thus
\bq \nb_x div(f)=0.
\label{a_div}
\eq 
Substituting (\ref{a_f}) and (\ref{a_div}) into (\ref{a_cond2}), we get
\[\frac{\dot{\al}+\dot{\be}}{\al+\be}\nb_x\log p +\frac{\al\dot{\be}-\dot{\al}\be}{\al+\be}\nb_x\log h=-(\nb_x\nb_x^T\log p)[K_1\nb_x\log p+K_2\nb_x\log h] \]
\bq
-[(\nb_x\nb_x^T\log p)K_1^T+(\nb_x\nb_x^T\log h) K_2^T](\nb_x\log p)+(\nb_x\nb_x^T\log p)Q(\nb_x\log p).
\label{a_cond3}
\eq
By setting the coefficient matrices corresponding to $\nb_x\log p$ and $\nb_x\log h$ equal on both sides of (\ref{a_cond3}), we obtain
\bq
\frac{\al\dot{\be}-\dot{\al}\be}{\al+\be}I = -(\nb_x\nb_x^T\log p)K_2,
\label{a_eq_K1}
\eq
and
\bq
\frac{\dot{\al}+\dot{\be}}{\al+\be}I = -(\nb_x\nb_x^T\log p)K_1
-(\nb_x\nb_x^T\log p)K_1^T-(\nb_x\nb_x^T\log h) K_2^T+(\nb_x\nb_x^T\log p)Q.
\label{a_eq_K2}
\eq
Under the assumption (A2), $\nb_x\nb_x^T\log p$ is invertible. It follows from (\ref{a_eq_K1}) that
\bq  
K_2 = - \frac{\al\dot{\be}-\dot{\al} \be}{\al+\be}(\nb_x\nb_x^T\log p)^{-1}.
\label{a_K2}
\eq
Substituting this form of $K_2$ into (\ref{a_eq_K2}), and also re-arranging terms, we have
\[
(\nb_x\nb_x^T\log p)(K_1+K_1^T)=(\nb_x\nb_x^T\log p)Q
+\frac{\al\dot{\be}-\dot{\al} \be}{\al+\be}(\nb_x\nb_x^T\log h) (\nb_x\nb_x^T\log p)^{-1}
\]
\[
-\frac{\dot{\al}+\dot{\be}}{\al+\be}I.
\]
Solving for $K_1+K_1^T$ gives
\[
K_1+K_1^T=Q
+\frac{\al\dot{\be}-\dot{\al} \be}{\al+\be}(\nb_x\nb_x^T\log p)^{-1}(\nb_x\nb_x^T\log h) (\nb_x\nb_x^T\log p)^{-1}
\]
\[
-\frac{\dot{\al}+\dot{\be}}{\al+\be}(\nb_x\nb_x^T\log p)^{-1}.
\]
Setting $K_1$ to be symmetric, we know from the previous equations that $K_1$ is given by
\[
K_1=\frac{1}{2}Q
+\frac{\al\dot{\be}-\dot{\al} \be}{2(\al+\be)}(\nb_x\nb_x^T\log p)^{-1}(\nb_x\nb_x^T\log h) (\nb_x\nb_x^T\log p)^{-1}
\]
\bq
-\frac{\dot{\al}+\dot{\be}}{2(\al+\be)}(\nb_x\nb_x^T\log p)^{-1}.
\label{a_K1}
\eq
The combination of (\ref{a_K1}) and (\ref{a_K2}) gives (\ref{f1})-(\ref{K2}). Q.E.D.
\vs

\noindent
{\sc \bf Proof of Theorem 2.2}

Under the assumption (A1), $\nb_x\log p_0$ and $\nb_x\log p_1$ are linear in $x$. We know from (\ref{a_logp1}) that $\nb_x\log h$ is also linear in $x$. Therefore, we may write 
\bq
\log p_0 = \frac{1}{2}x^TA_0x+b_0^Tx+c_0,
\label{logp0}
\eq
%and
\bq
\log h = \frac{1}{2} x^TA_hx+b_h^Tx+c_h,
\label{logh}
\eq
where $A_0, A_h\in \mathbb{R}^{n\times n}, b_0, b_h\in \mathbb{R}^n, c_0, c_h\in \mathbb{R}$ are constant matrices, vectors, and scalars, respectively. They are not functions of either $x$ or $\lambda$. Note that $A_0$ and $A_h$ are symmetric because they are Hessian matrices of density functions.
Consequently, $\log p$ has the form
\[
\log p = (\al+\be)\log p_0+\be\log h-\log\Gamma(\lambda)-\be\log \Gamma(1)
\]
\[= \frac{1}{2}x^T[(\al+\be)A_0+\be A_h]x+[(\al+\be)b_0+\be b_h]^Tx+(\al+\be)c_0+\be c_h-\log\Gamma(\lambda)-\be\log \Gamma(1).
\]
The gradient and Hessian of $\log h$ and $\log p$ are, respectively
\[
\nb_x\log h = A_hx+b_h, 
\nb_x\nb_x^T\log h = A_h,
\]
\[
\nb_x\log p = [(\al+\be)A_0+\be A_h]x+(\al+\be)b_0+\be b_h,
\nb_x\nb_x^T\log p = (\al+\be)A_0+\be A_h.
\]
The assumption (A2) says that $(\al+\be)A_0+\be A_h$ is non-singular for all $\lambda\in [0, 1]$. In this case, for each $\lambda\in[0, 1]$, $p(x,\lambda)$ in (\ref{loghom}) is Gaussian distributed with mean $x_\mu(\lambda)$ determined by setting $\nb_x\log p|_{x=x_{\mu}} = 0$, which gives
\bq
x_\mu(\lambda)=-[(\al+\be)A_0+\be A_h]^{-1}[(\al+\be)b_0+\be b_h],
\label{xmu}
\eq
and its covariance matrix $P_\mu(\lambda)$ is given by
\bq
P_\mu(\lambda) = -(\nb_x\nb_x^T\log p)^{-1} = - [(\al+\be)A_0+\be A_h]^{-1}.
\label{Pmu}
\eq
On the other hand, under the assumption (A1), the flow equation (\ref{flow}) is a linear stochastic process
\bq
dx=(Fx+b)d\lambda + qdw_{\lambda}
\label{linear}
\eq
in which
\[ 
F=\nb_x f = K_1\nb_x\nb_x^T\log p+K_2\nb_x\nb_x^T\log h, \hs
b=f-Fx.
\]
For the linear stochastic differential equation (\ref{linear}), the mean $\bar{x}$ and covariance matrix $P$ of $x$ satisfy, respectively, \cite{Arn}
\bq
\frac{d \bar{x}}{d\lambda}=F\bar{x}+b,
\label{mean}
\eq
\bq
\frac{d P}{d\lambda} = FP+PF^T+Q.
\label{cov}
\eq
Note that $x(0)$ is Gaussian. We know from the theory of linear stochastic differential equations that the probability density function of $x$  (\ref{linear}) is Gaussian for each $\lambda$ \cite{Arn}. To prove that the flow (\ref{flow}), or equivalently (\ref{linear}), has desired density function $p(x,\lambda)$ (\ref{loghom}), we only need to verify that $x_\mu$ and $P_\mu$ defined in (\ref{xmu})-(\ref{Pmu}) indeed satisfy (\ref{mean}) and (\ref{cov}), respectively.
\vs

We first consider the mean $x_\mu$. For notational convenience, denote $S(\lambda)=\nb_x\nb_x^T\log p$. Then
\[
 \frac{dS}{d\lambda}=(\dot{\al}+\dot{\be})A_0+\dot{\be} A_h.
\]
For $\bar{x}=x_\mu$, using Lemma A.1,
%using the equality for differentiating an invertible matrix $B(\theta)$,
%\[
%dB^{-1}/d\theta=-B^{-1}(dB/d\theta)B^{-1},
%\]  
the left hand side (LHS) of (\ref{mean}) is
\[
\textrm{LHS of (\ref{mean})}=\frac{d x_{\mu}}{d\lambda} =\frac{d}{d\lambda}\{-S^{-1}[(\al+\be)b_0+\be b_h]\}
\]
\[
=-(\frac{dS^{-1}}{d\lambda})[(\al+\be)b_0+\be b_h]-S^{-1}[(\dot{\al}+\dot{\be})b_0+\dot{\be} b_h]
\]
\[
=S^{-1}[(\dot{\al}+\dot{\be})A_0+\dot{\be} A_h)]S^{-1}[(\al+\be)b_0+\be b_h]-S^{-1}[(\dot{\al}+\dot{\be})b_0+\dot{\be} b_h]
\]
\bq
=-S^{-1}[(\dot{\al}+\dot{\be})A_0+\dot{\be} A_h)]x_{\mu}-S^{-1}[(\dot{\al}+\dot{\be})b_0+\dot{\be} b_h].
\label{mean_LHS}
\eq
On the other hand, (\ref{a_d2logpdlam_1}) and (\ref{a_d2logpdlam}) are two expressions of $\nb_x[\frac{\partial \log p}{\partial \lambda}]$. Therefore,
\[
\frac{\dot{\al}+\dot{\be}}{\al+\be}\nb_x\log p +\frac{\al\dot{\be}-\dot{\al} \be}{\al+\be}\nb_x\log h= (\dot{\al}+\dot{\be})\nb_x\log p_0 +\dot{\be}\nb_x\log h
\]
or by re-arranging the terms
\bq
\frac{\al\dot{\be}-\dot{\al}\be}{\al+\be}\nb_x\log h= 
-\frac{\dot{\al}+\dot{\be}}{\al+\be}\nb_x\log p+
(\dot{\al}+\dot{\be})\nb_x\log p_0 +\dot{\be}\nb_x\log h.
\label{a_rhs1}
\eq
Note that $x_\mu$ satisfies $\nb_x\log p |_{x=x_\mu} = 0$ by construction. Setting $\bar{x}=x_\mu$, also using the relationship (\ref{a_rhs1}), the right hand side (RHS) of (\ref{mean}) becomes
\[
\textrm{RHS of (\ref{mean})} = f |_{\bar{x}=x_\mu}
=S^{-1}(K_1\nb_x\log p+K_2\nb_x\log p)|_{\bar{x}=x_\mu}
\]
\[
=- \frac{\al\dot{\be}-\dot{\al}\be}{\al+\be}S^{-1}(\nb_x\log h)|_{\bar{x}=x_{\mu}}
\]
\[
=-S^{-1}[(\dot{\al}+\dot{\be})\nb_x\log p_0+\dot{\be}\nb_x\log h]|_{\bar{x}=x_{\mu}}
\]
\[
=-S^{-1}[(\dot{\al}+\dot{\be})(A_0x_{\mu}+b_0)+\dot{\be}(A_hx_{\mu}+b_h)]
\]
\[
=-S^{-1} [(\dot{\al}+\dot{\be})A_0+\dot{\be} A_h)]x_{\mu}-S^{-1}[(\dot{\al}+\dot{\be})b_0+\dot{\be} b_h] 
\]
\[
= \textrm{LHS of (\ref{mean})}
\]
according to (\ref{mean_LHS}).

For the covariance matrix $P_\mu$, we have $P_\mu=-S^{-1}$. Therefore,
\bq
\textrm{LHS of (\ref{cov})} = \frac{d P_{\mu}}{d\lambda}=S^{-1} [(\dot{\al}+\dot{\be})A_0+\dot{\be} A_h)] S^{-1}
\label{cov_LHS}
\eq
and
\[
\textrm{RHS of (\ref{cov})} =FP_\mu+P_\mu F^T+Q
=[K_1S+K_2 A_h]P_\mu+P_\mu [S K_1^T+A_h K_2^T]+Q
\]
\[
=\frac{\al\dot{\be}-\dot{\al} \be}{\al+\be}S^{-1}A_hS^{-1}+\frac{\dot{\al}+\dot{\be}}{\al+\be}S^{-1}
=S^{-1}[\frac{\al\dot{\be}-\dot{\al} \be}{\al+\be}A_h+\frac{\dot{\al}+\dot{\be}}{\al+\be}S]S^{-1}
\]
\[
=S^{-1}[(\dot{\al}+\dot{\be})A_0+\dot{\be} A_h]S^{-1}=\textrm{RHS of (\ref{cov})}.
\]
Q.E.D.
\vs

\noindent
{\sc \bf Proof of Lemma 2.2}

Taking derivative with respect to $\lambda$ on both sides of (\ref{a_ddlogp}), 
\[
\frac{d}{d\lambda}\nb_x\nb_x^T\log p=(\dot{\al}+\dot{\be})\nb_x\nb_x^T\log p_0+\dot{\be}\nb_x\nb_x^T\log h.
\]
Therefore, for $M=-\nb_x\nb_x^T\log p$, %we know from the previous equation that
\bq
\frac{dM}{d\lambda} = -(\dot{\al}+\dot{\be})\nb_x\nb_x^T\log p_0-\dot{\be}\nb_x\nb_x^T\log h. 
\label{a_dM1}
\eq
By taking the gradient again on both sides of (\ref{a_dlogp0}), we obtain
\bq
\nb_x\nb_x^T\log p_0 = \frac{1}{\al+\be}\nb_x\nb_x^T\log p-\frac{\be}{\al+\be}\nb_x\nb_x^T\log h.
\label{a_dM2}
\eq
Combining (\ref{a_dM1}) and (\ref{a_dM2}) gives
\bq
\frac{dM}{d\lambda} = -\frac{\dot{\al}+\dot{\be}}{\al+\be}\nb_x\nb_x^T\log p-\frac{\al\dot{\be}-\dot{\al}\be}{\al+\be}\nb_x\nb_x^T\log h. 
\label{a_dM}
\eq
According to (\ref{linearized}), 
\[
dV= \tilde{x}^T[F^TM+MF+\frac{dM}{d\lambda}]\tilde{x}d\lambda.
\]
Substituting the form of $F$ (\ref{Jacobian1}) and $dM/d\lambda$ (\ref{a_dM}) into the previous equation, 
\[
dV= \tilde{x}^T[\frac{1}{2}Q(\nb_x\nb_x^T\log p)
-\frac{\al\dot{\be}-\dot{\al} \be}{2(\al+\be)}(\nb_x\nb_x^T\log p)^{-1}(\nb_x\nb_x^T\log h) -\frac{\dot{\al}+\dot{\be}}{2(\al+\be)}I]^TM\tilde{x}d\lambda
\]
\[+\tilde{x}^TM[\frac{1}{2}Q(\nb_x\nb_x^T\log p)
-\frac{\al\dot{\be}-\dot{\al} \be}{2(\al+\be)}(\nb_x\nb_x^T\log p)^{-1}(\nb_x\nb_x^T\log h) -\frac{\dot{\al}+\dot{\be}}{2(\al+\be)}I]\tilde{x} d\lambda
\]
\[
+\tilde{x}^T[-\frac{\dot{\al}+\dot{\be}}{\al+\be}\nb_x\nb_x^T\log p-\frac{\al\dot{\be}-\dot{\al}\be}{\al+\be}\nb_x\nb_x^T\log h]\tilde{x}d\lambda
\]
\[
= \tilde{x}^T\{[-\frac{1}{2}QM +\frac{\al\dot{\be}-\dot{\al} \be}{2(\al+\be)}M^{-1}(\nb_x\nb_x^T\log h) -\frac{\dot{\al}+\dot{\be}}{2(\al+\be)}I]^TM
+M[-\frac{1}{2}QM
\]
\[
+\frac{\al\dot{\be}-\dot{\al} \be}{2(\al+\be)}M^{-1}(\nb_x\nb_x^T\log h) -\frac{\dot{\al}+\dot{\be}}{2(\al+\be)}I]
+[\frac{\dot{\al}+\dot{\be}}{\al+\be}M-\frac{\al\dot{\be}-\dot{\al}\be}{\al+\be}\nb_x\nb_x^T\log h]\}\tilde{x}d\lambda
\]
\[
=-\tilde{x}^TMQM\tilde{x}d\lambda
\]
which is (\ref{dV}). Q.E.D.
\vs

\noindent
{\sc \bf Proof of Theorem 2.3} 

Note that $M(\lambda)$ is symmetric and 
\[
M(0)=-\al(0)\nb_x\nb_x^T\log p_0-\be(0)\nb_x\nb_x^T\log p_1=-\nb_x\nb_x^T\log p_0
\]
because $\al(0)=1$ and $\be(0)=0$ according to (\ref{albe_bc}). For any $Q\geq 0$, $MQM\geq 0$. We know from (\ref{dV}) that $dV\leq 0$ for all $\lambda$. Consequently, 
\[
V\leq V|_{\lambda=0}=\tilde{x}_0^TM(0)\tilde{x}_0=
\tilde{x}_0^T(-\nb_x\nb_x^T\log p_0)\tilde{x}_0
=c.
\] 
On the other hand, under the assumption (A3), $M\geq M_0$. Therefore, 
\[
\tilde{x}^TM_0\tilde{x}\leq V\leq c,
\]
which is (1). 

To prove (2), note that $Q>0$ leads to
\bq
\tilde{x}^T(MQM)\tilde{x}\geq \lambda_{min}(Q) \tilde{x}^TMM\tilde{x}.
\label{MQM}
\eq
Under the assumption (A3), $M$ is positive definite. We can write $M=S_1^2$ with $S_1=M^{1/2}$ positive definite. Then 
\[
\tilde{x}^T MM\tilde{x}=\tilde{x}^TS_1 S_1^2 S_1\tilde{x}=\tilde{x}^TS_1M S_1\tilde{x}\geq\tilde{x}^TS_1M_0S_1\tilde{x}\geq \lambda_{min}(M_0)\tilde{x}^TS_1S_1\tilde{x}
\]
\[
=\lambda_{min}(M_0)\tilde{x}^TM\tilde{x}=\lambda_{min}(M_0)V.
\]
Combining the previous inequality with (\ref{dV}) and (\ref{MQM}) leads to
\[ 
dV =-\tilde{x}^T(MQM)\tilde{x}d\lambda \leq -\lambda_{min}(Q) \tilde{x}^TMM\tilde{x}d\lambda \leq -\lambda_{min}(Q)\lambda_{min}(M_0)Vd\lambda =-rV d\lambda.
\]
Applying the Gronwall-Bellman Inequality in Lemma A.2 to the previous inequality gives
\[
V \leq (V|_{\lambda=0}) e^{-r\lambda}=c e^{-r\lambda}, \forall \lambda \in [0, 1].
\]
On the other hand, under the assumption (A3), $V=\tilde{x}^TM\tilde{x}\geq \tilde{x}^TM_0\tilde{x}$, which, together with the previous inequality, gives (\ref{exp_bound}). 
Q.E.D.
\vs

\noindent
{\sc \bf Proof of Theorem 3.1}

We obtain the optimal solution of $u$ by applying Pontryagin maximum principle \cite{BH}.
For the optimal control problem (\ref{state})-(\ref{obj}), the corresponding Hamiltonian function is
\[
H(\be,u,\xi) = \frac{1}{2}u^2+\mu \kappa_{\nu}(M)+\xi u
\]
in which $\xi$ is the adjoint variable (costate). Let $u^*(\lambda)$, $\be^*(\lambda)$ and $\xi^*$ be the optimal solutions of control, state and adjoint variable, respectively. The Pontryagin maximum principle states that the follow three necessary conditions are satisfied \cite{BH}.

\noindent 1. The optimal control $u^*$ satisfies
\[
H(\be^*,u^*,\xi^*) = \min_{u} H(\be^*,u,\xi^*).
\]
Setting the partial derivative of $H(\be^*,u,\xi^*)$ with respect to $u$ to 0, we obtain
\[
0=\frac{\partial H}{\partial u} = u+\xi^*
\]
which gives the optimal control
\bq
u^* = - \xi^*.
\label{opti_u}
\eq
2. The adjoint variable satisfies the adjoint equation
\bq
\frac{d\xi^*}{d\lambda}=-\frac{\partial H(\be^*,u^*,\xi^*)}{\partial \be}=-\mu\frac{\partial \kappa_{\nu}(M)}{\partial \be}|_{\be=\be^*}.
\label{opti_xi}
\eq
3. The optimal state $\be^*$ satisfies
\bq
\frac{d\be^*}{d\lambda}=u^*(\lambda)=-\xi^*, 
\label{opti_b1}
\eq
\[
\be^*(0)=0, \be^*(1)=1.
\]
Combining (\ref{opti_xi}) and (\ref{opti_b1}) gives (\ref{opti_b}). The optimal solution $\be^*$ must satisfy the boundary condition (\ref{opti_bv}). Q.E.D.
\vs

\noindent
{\sc \bf Proof of Theorem 3.2} 

We prove by contradiction. Let $\be^*$ be a continuously differentiable solution to the optimal control (\ref{state})-(\ref{obj}) with corresponding $u^*=\dot{\be}^*$. The boundary condition (\ref{BV}) gives that $\be^*(0)=0, \be^*(1)=1$. Assume that there exists a $\lambda^+\in (0, 1)$ such that $\be^*(\lambda^+)<0$. % which we shall prove impossible by contradiction.
Define
\[
\lambda^+_1=\max\{ \lambda | \be^*(\lambda)=0, \lambda\in [0,\lambda^+)\}.
\]
Note that $\lambda^+_1$ is non-empty (and thus well-defined) because $\be^*(0)=0$ and $\lambda^+_1<\lambda^+$. Since $\be^*$ is continuous over $[0,1 ]$, we must have
\bq
\be^*(\lambda) < 0, \forall \lambda \in (\lambda^+_1, \lambda^+].
\label{a_beta_1}
\eq
By assumption, $-\nb_x\nb_x^T\log p_0>0$. The continuity of $\be^*$ and $\be^*(\lambda^+_1)=0$ assures that there exists a $\lambda^+_2$ sufficiently close to $\lambda^+_1$ and $\lambda^+_2 \in (\lambda^+_1,\lambda^+]$ (thus $\lambda^+_2>\lambda^+_1$) such that
\bq
-\nb_x\nb_x^T\log p_0+\be^*(-\nb_x\nb_x^T\log h)>0, \forall \lambda \in (\lambda^+_1,\lambda^+_2).
\label{a_beta_2}
\eq
Because $\be^*$ is continuously differentiable, $u^*=\dot{\be}^*$ is continuous over $(\lambda^+_1,\lambda^+_2)$. We must have 
\bq
\int_{\lambda^+_1}^{\lambda^+_2} \frac{1}{2}(u^*)^2d\lambda>0,
\label{a_u}
\eq
%$u^*\not\equiv 0, \forall\lambda \in (\lambda^+_1,\lambda^+_2)$. 
otherwise we would have $u^*\equiv 0$ and consequently $\be^*\equiv 0, \forall\lambda \in(\lambda^+_1,\lambda^+_2)$ because $\be^*(\lambda^+_1)=0$. 
%We also have
%\bq
%-\be^*(\lambda)>0,  \forall \lambda\in (\lambda^+_1,\lambda^+_2).
%\label{a_beta_1new}
%\eq
%according to (\ref{a_beta_1}).

We know from the assumption $\kappa_{\nu}(\nb_x\nb_x^T\log h)\leq\kappa_{\nu}(\nb_x\nb_x^T\log p_0)$ that
$\kappa_\nu(-\nb_x\nb_x^T\log h)\leq \kappa_\nu(-\nb_x\nb_x^T\log p_0)$.
For a monotone norm $||\cdot||_\nu$ and
for each $\lambda \in (\lambda^+_1,\lambda^+_2)$, by applying Lemma A.3 with $d_1=0, d_2=-\be^*(\lambda)>0, A=-\nb_x\nb_x^T\log p_0>0, B=-\nb_x\nb_x^T\log h>0$, we know from (\ref{a_beta_2}) and (\ref{a_u}) that
\bq
\kappa_\nu(-\nb_x\nb_x^T\log p_0)
\leq 
\kappa_\nu(-\nb_x\nb_x^T\log p_0+\be^*(-\nb_x\nb_x^T\log h)) ,
\forall \lambda \in (\lambda^+_1,\lambda^+_2).
\label{a_beta_3}
\eq

Define
\[
\be^+(\lambda) = 
\left\{
\ba{ll}
0, & \textrm{if $\lambda\in (\lambda^+_1,\lambda^+_2)$}, \\
\be^*(\lambda), &  \textrm{if $\lambda\notin (\lambda^+_1,\lambda^+_2)$}.
\ea
\right.
\] 
Then $\be^+(\lambda)$ and $\be^*(\lambda)$ are not identical over $[0, 1]$:
\[
\be^+(\lambda)=0>\be^*(\lambda), \forall \lambda \in (\lambda^+_1,\lambda^+_2). 
\]
For this $\be^+(\lambda)$, denote $u^+=d\be^+/d\lambda$. According to (\ref{a_u})
\bq
\int_{\lambda^+_1}^{\lambda^+_2} \frac{1}{2}(u^+)^2d\lambda
%= \int_{\lambda^+_1}^{\lambda^+_2} \frac{1}{2}(\frac{d\be^+}{d\lambda})^2d\lambda 
=0
< \int_{\lambda^+_1}^{\lambda^+_2} \frac{1}{2}(u^*)^2d\lambda.
\label{a_u1}
\eq
On the other hand, let $M^+=-\nb_x\nb_x^T\log p_0+\be^+(-\nb_x\nb_x^T\log h)$ and $M^*=-\nb_x\nb_x^T\log p_0+\be^*(-\nb_x\nb_x^T\log h)$. An equivalent form of (\ref{a_beta_3}) is
\bq
\kappa_\nu(M^+) \leq \kappa_\nu(M^*), \forall \lambda \in (\lambda^+_1,\lambda^+_2).
\label{a_kappa1}
\eq
The combination of (\ref{a_u1}) and (\ref{a_kappa1}) gives
\bq
\int_{\lambda^+_1}^{\lambda^+_2} [\frac{1}{2}(u^+)^2+\mu \kappa_\nu(M^+)]d\lambda 
< \int_{\lambda^+_1}^{\lambda^+_2} [\frac{1}{2}(u^*)^2+\mu \kappa_\nu(M^*)]d\lambda.
\label{a_inequality1}
\eq
Consequently, also noticing that $\be^+\equiv\be^*$ and $(u^+)^2 \equiv (u^*)^2$ for $\lambda\in [0, \lambda^+_1]\cup [\lambda^+_2, 1]$,
\[
J(\be^+,u^+) 
\]
\[
=\int_0^1 [\frac{1}{2}(u^+)^2+\mu \kappa_\nu(M^+)]d\lambda
\]
\[
=\int_0^{\lambda^+_1} [\frac{1}{2}(u^+)^2+\mu \kappa_\nu(M^+)]d\lambda 
+\int_{\lambda^+_1}^{\lambda^+_2} [\frac{1}{2}(u^+)^2+\mu \kappa_\nu(M^+)]d\lambda 
+\int_{\lambda^+_2}^1 [\frac{1}{2}(u^+)^2+\mu \kappa_\nu(M^+)]d\lambda 
\]
\[
=\int_0^{\lambda^+_1} [\frac{1}{2}(u^*)^2+\mu \kappa_\nu(M^*)]d\lambda 
+\int_{\lambda^+_1}^{\lambda^+_2} [\frac{1}{2}(u^+)^2+\mu \kappa_\nu(M^+)]d\lambda 
+\int_{\lambda^+_2}^1 [\frac{1}{2}(u^*)^2+\mu \kappa_\nu(M^*)]d\lambda 
\]
\[
< \int_0^{\lambda^+_1} [\frac{1}{2}(u^*)^2+\mu \kappa_\nu(M^*)]d\lambda 
+\int_{\lambda^+_1}^{\lambda^+_2} [\frac{1}{2}(u^*)^2+\mu \kappa_\nu(M^*)]d\lambda 
+\int_{\lambda^+_2}^1 [\frac{1}{2}(u^*)^2+\mu \kappa_\nu(M^*)]d\lambda 
\]
\[
= \int_0^{1} [\frac{1}{2}(u^*)^2+\mu \kappa_\nu(M^*)]d\lambda 
\]
\bq
= J(\be^*,u^*),
\label{a_inequality2}
\eq
%The value of the objective function for $\be^+$ is less that for $\be^*$,
which contradicts the assumption that $\be^*$ minimizes $J$. This contradiction means that there does not exist a $\lambda^+$ such that $\be^*(\lambda^+)<0$. In other words, we must have $\be^*(\lambda)\geq 0$ for all $\lambda \in [0, 1]$.  
Q.E.D.

%\newpage

\end{document}